\documentclass[12pt]{article}
\usepackage{amsmath,amsthm,verbatim,amssymb,amsfonts,amscd,diagrams, graphicx, mathrsfs, hyperref,mathtools,tipa,multicol, enumitem}
\usepackage[square,sort,comma,numbers]{natbib}
\usepackage[usenames,dvipsnames]{color}
\usepackage{amstext}
\usepackage{color}

\theoremstyle{plain}
\newtheorem{theorem}{Theorem}
\newtheorem*{theorem*}{Theorem}

\theoremstyle{definition}

\theoremstyle{remark}
\newtheorem{remark}[theorem]{Remark}

\newcommand{\codim}{\text{codim}}
\newarrow{ul}---->
\newarrow{Backwards}<----

\newcommand{\into}{\hookrightarrow}
\newcommand{\Z}{\mathbb{Z}}

\newcommand{\R}{\mathbb{R}}

\newcommand{\PP}{\mathbb{P}}

\renewcommand{\H}{\mathbb H}

\newcommand{\mc}[1]{\mathcal{#1}}

\newcommand{\mf}{\mathfrak}

\newarrow{Onto}----{>>}
\newarrow{Equals}=====

\newcommand{\xr}{\xrightarrow}

\begin{document}

\title{Two short proofs of the topological invariance of intersection homology}
\author{Greg Friedman}

\date{May 22, 2019}

\maketitle


\begin{abstract}
We indicate two  short proofs of the Goresky-MacPherson topological invariance of intersection homology. One proof is very short but requires the Goresky-MacPherson support and cosupport axioms; the other is slightly longer but does not require these axioms and so is adaptable to more general perversities.  {\let\thefootnote\relax\footnote{{

}}} 
\end{abstract}

While the key feature of the intersection homology of Goresky and MacPherson \cite{GM1,GM2} is undoubtedly that it can be used to extend Poincar\'e duality to singular spaces, a second important property  is that it is a topological invariant, at least if one restricts to the original perversity parameters introduced in \cite{GM1}. This means that while the definition of the intersection homology groups depends on a choice of topological stratification, the resulting groups do not. 

This invariance was originally proven by Goresky and MacPherson in \cite{GM2} using the sheaf-theoretic approach to intersection homology on stratified pseudomanifolds. They showed that the Deligne sheaf complexes, whose hypercohomology gives intersection homology, are characterized up to quasi-isomorphism by a simple set of axioms  and then moved through various equivalent sets of axioms, ultimately finding ones that do not depend on the specific stratification but essentially only on the constructibility of the sheaf complex with respect to \emph{some} stratification  and conditions on dimensions of supports and ``cosupports'' of the derived cohomology sheaves. The crux of the argument is then the construction of a single sheaf complex that satisfies all of these axioms and so is quasi-isomorphic to those constructed for any specific stratification. As noted in \cite[Section V.4]{Bo}, where a more detailed treatment is given by Borel, the main difficulty is the construction of an appropriate stratification that is then used to construct this universal example. 

Alternative and simpler proofs followed. One key simplification that addresses Borel's concern is to replace the stratification used by Goresky-MacPherson, which is defined by sheaf-theoretic constructibility properties, with a more concrete topological stratification, namely the intrinsic stratification that a pseudomanifold possesses when considered in the more general class of Siebenmann's CS sets \cite{Si72}. The CS sets possess quite natural intrinsic stratifications determined by an equivalence relation so that $x\sim y$  if  there are homeomorphic neighborhood pairs $(U,x)\cong (V,y)$; see  \cite[Section 2.10]{GBF35}. The idea of recasting intersection homology and its topological invariance in the class of CS sets is due to King \cite{Ki}, who provided a sheaf-less proof of the  topological invariance of singular chain intersection homology in that context. However, King's proof is also somewhat intricate, involving an intertwined induction argument on three different statements. 

In \cite{HS91}, Habegger and Saper extended sheaf-theoretic intersection homology to CS sets, generalizing to \emph{codimension $\geq c$ intersection cohomology theories}, \emph{$c$-ICTs} for short (see \cite[Section 11]{GBF26} for the relation between these theories and other generalizations of intersection homology). The perspective here is a bit different, with the topological invariance of $c$-ICTs being more or less built in as an axiom (see \cite[Remark 4.2 and the proof of Theorem 6.2]{HS91}). It is then shown that the collection of sheaves satisfying the support and cosupport conditions of Goresky-MacPherson constitutes a $c$-ICT \cite[Definition 4.3 and Propostion 4.5]{HS91}. From this it is possible to deduce the topological invariance of Goresky and MacPherson's intersection homology, essentially by the observation at the core of our first proof below, though we will show how to make a direct argument without invoking the additional machinery of \cite{HS91}.

More recently, Chataur, Saralegi-Aranguren, and Tanre \cite{CST-inv} have given another proof of topological invariance for much more general perversities, once again without sheaves. Their proof utilizes a Mayer-Vietoris argument together with some of King's topological arguments, though in the case of Goresky-MacPherson perversities their proof considerably simplifies King's. 

Our goal here is to provide two fairly direct and very short sheaf-theoretic proofs of the topological invariance of intersection homology in the classical setting of Goresky-MacPherson perversities and constant coefficients, building on some of the foundation of the Goresky-MacPherson proof. Our first proof is inspired by \cite{HS91}, though it is more direct since we do not require the $c$-ICT machinery. This proof does involve the support and cosupport conditions of Goresky and MacPherson, which means that it does not readily extend to more exotic perversities than the Goresky-MacPherson perversities (see \cite{GBF43} for further discussion). The second proof requires slightly more work but readily generalizes not only to more general perversities but to the ``torsion-sensitive intersection homology'' introduced in \cite{GBF32}. See Remark \ref{R: general}, below, for an indication of handling more general perversities  and \cite{GBF43} for the torsion-sensitive case.

To keep this note as brief as possible, we focus almost entirely on constant coefficients and Goresky-MacPherson perversities, i.e.\ the original theorem of Goresky-MacPherson \cite[Theorem 4.1]{GM2}, though see Remark \ref{R: twisted}.
We also refer to other sources for more detailed background, especially \cite{GBF35} for CS sets, Borel \cite{Bo} for sheaf-theoretic intersection homology and the axiomatic approach to it, and Sch\"urmann \cite[Chapter 4]{Sch03} for the preservation of constructibility by pushforwards and pullbacks (with or without compact support) in the context of CS sets (though see also Habbeger-Saper \cite[Appendix]{HS91} for constructibility on CS sets).

\bigskip

To fix notation: Let $X$ be a paracompact CS set. Such spaces are  Hausdorff (by definition), locally compact \cite[Lemma 2.3.15]{GBF35}, metrizable \cite[Proposition 1.11]{CST-inv}, and of finite cohomological dimension (\cite[Lemma 6.3.46]{GBF35} and \cite[Theorem II.16.8]{Br}). We assume $X$ filtered by $X=X^n\supset X^{n-2}\supset X^{n-3}\supset\cdots\supset X^{-1}=\emptyset$ and that $X-X^{n-2}$ is dense. Let $U_k=X-X^{n-k}$, let $i_k:U_k\into U_{k+1}$, and let $S_k=X^k-X^{k-1}$ be the union of $k$-dimensional \emph{strata}.
 Let $R$ be a commutative Noetherian ring of finite cohomological dimension, and let $\mc E$ be a constant sheaf on $X$ whose stalks are finitely-generated $R$-modules. Let $\bar p$ be a GM-perversity, i.e.\ $\bar p(2)=0$ and $\bar p(k)\leq \bar p(k+1)\leq \bar p(k)+1$. In this case the Goresky-MacPherson-Deligne sheaf is defined to be 
$\mc P^*=    \tau_{\leq \bar p(n)}Ri_{n*}      \cdots \tau_{\leq \bar p(2)}Ri_{2*}\mc E_{U_2}.$
Let $\mf X$ denote the intrinsic CS set stratification of $X$, which coarsens $X$ (every stratum of $X$ is contained in a stratum of $\mf X$), and let $\PP^*$ be the analogous $(\bar p,\mc E)$-Deligne sheaf with respect to $\mf X$.
For sheaf complexes, the symbol $\cong$ denotes isomorphism in the derived category, i.e.\ quasi-isomorphism. Recall that a sheaf complex is called $X$-clc (for \emph{cohomologically locally constant}) if its derived cohomology sheaves are locally constant on each stratum. On a CS set, if $j$ is any inclusion of a locally closed subset that is a union of strata then $j^*$, $j^!$, $j_!$, and $Rj_*$ all preserve this property of \emph{constructibility} by \cite[Proposition 4.0.2.3]{Sch03} (see also \cite[Proposition 4.2.1.2.b]{Sch03}).

We provide two proofs for the following theorem. 

\begin{theorem*}\label{T: theorem}
$\mc P^*$ is quasi-isomorphic to $\PP^*$. Consequently the Deligne sheaves with respect to any two CS set stratifications are quasi-isomorphic. 
\end{theorem*}

The basic idea of both proofs involves a simplified application of the Goresky-MacPherson axiomatics: It is shown by Goresky-MacPherson \cite{GM2} (cf. Borel \cite[Section V]{Bo}) that there are various equivalent sets of axioms that characterize the Deligne sheaf with perversity $\bar p$ and coefficients $\mc E$ on an appropriately stratified space $X$.
Even though we will not utilize all of the axioms directly, we recall them briefly for the reader's benefit. As we have fixed $\mc E$ and $\bar p$, we omit them from the notation. In the statements, $\mc S^*$ is a complex of sheaves on $X$, we let $f_x:x\into X$ be the inclusion, $\bar q$ is the complementary perversity to $\bar p$ (i.e.\ $\bar q(k)=k-2-\bar p(k)$), and $\bar p^{-1}(i)=\min\{c\mid \bar p(c)\geq i\}$ (taking $\bar p^{-1}(i)=\infty$ if $i>\bar p(n)$). 

\begin{description}
\item[Ax1$(X)$:]\hfill

\begin{enumerate}

\item\label{I: normal} $\mc S^*$ is bounded, $\mc S^j=0$ for $j<0$, $\mc S^*|_{U_2}\cong \mc E|_{U_2}$,

\item\label{I: truncate} If $x\in S_{n-k}$, $k\geq 2$, then $H^j(\mc S_x)=0$ of $j>\bar p(k)$,

\item\label{I: attach} The attachment map $\mc S^*|_{U_{k+1}}\to Ri_{k*}\mc S^*|_{U_k}$ is a quasi-isomorphism up to degree $\bar p(k)$. 

\end{enumerate}

\item[Ax1'$(X)$:] Same as Ax1$(X)$ but adding that  $\mc S^*$ is $X$-clc to Axiom \ref{I: normal} and replacing axiom \ref{I: attach} with:

\setlist[description]{font=\normalfont\space}
\begin{description}
\item[3'.] If $x\in S_{n-k}$ then $H^j(f_x^!\mc S^*)=0$ for $j<n-\bar q(k)$.

\end{description}

\item[Ax2$(X)$:] \hfill
\begin{enumerate}

\item $\mc S^*$ is bounded, $\mc S^j=0$ for $j<0$, $\mc S^*|_{U_{2}}\cong \mc E|_{U_2}$, and $\mc S^*$ is $X$-clc,

\item $\dim\{x\in X\mid \mc H^j(f_x^*\mc S^*)\neq 0\}\leq n-\bar p^{-1}(i)$ for all $j>0$,

\item $\dim\{x\in X\mid \mc H^j(f_x^!\mc S^*)\neq 0\}\leq n-\bar q^{-1}(n-i)$ for all $j<n$.

\end{enumerate}

\end{description}

The axioms Ax1$(X)$ immediately follow from, and nearly immediately imply, a sheaf complex being the Deligne sheaf complex \cite[Theorem V.2.5]{Bo}. It is also elementary to show that 
if each stratum is a manifold, if the sheaf complex $\mc S^*$ is $X$-clc, and if $j^!\mc S^*$ is clc for each stratum inclusion $j:S\into X$, then Ax1$(X)$ is equivalent to  Ax1'$(X)$, which in turn is equivalent to  Ax2$(X)$ \cite[Proposition 4.3, Proposition 4.9]{Bo}. These results are proven in \cite{Bo}  under the assumption that the underlying space of $X$ is a pseudomanifold, but the arguments hold as well for our CS sets with dense $X-X^{n-2}$. Furthermore, all strata of CS sets are manifolds, and if $\mc S^*$ is $X$-clc then automatically $j^!\mc S^*$ is clc for each stratum inclusion $j:S\into X$ by \cite[Proposition 4.0.2.3]{Sch03}. Similarly, \cite[Proposition 4.0.2.3]{Sch03} implies that any sheaf complex satisfying Ax1$(X)$ is $X$-clc. So, for a given stratification of a CS set, a sheaf complex satisfies all three sets  of axioms  and is quasi-isomorphic to the Deligne sheaf if and only if it satisfies one of the sets of axioms.

The proofs of topological invariance in \cite{GM2,Bo} proceed by constructing a complex that satisfies Ax2 for \emph{every} pseudomanifold stratification of the underlying space. By the results cited in the preceding paragraph, this complex is then quasi-isomorphic to the Deligne sheaf with respect to each such stratification, and so they are all quasi-isomorphic to each other. This ``universal Deligne sheaf'' is just the Deligne sheaf with respect to a specially-constructed stratification, but, as noted above, the hard part is constructing that stratification. By contrast, the intrinsic CS set stratification is much simpler to construct; compare the constructions in \cite[Section V.4]{Bo} and \cite[Section 1]{Ki}. The key idea of our proofs is to show that in the context of CS sets the Deligne sheaf $\PP^*$ with respect to the intrinsic CS set stratification plays the same universal role, i.e.\ it is quasi-isomorphic to the Deligne sheaf with respect to any CS set stratification. 

With these preliminaries, our first proof of the Theorem is especially simple:

\begin{proof}[Proof 1]
By construction, the Deligne sheaf $\PP^*$ satisfies the axioms Ax1$(\mf X)$.  So by the above-cited results, $\PP^*$ also satisfies Ax2$(\mf X)$. But the only axiom in Ax2 that depends on the stratification is the first one, and as $\mf X$ coarsens $X$ we have that $\PP^*$ is also $X$-clc while $X-X^{n-2}\subset \mf X-\mf X^{n-2}$. Thus $\PP^*$ satisfies Ax2$(X)$. It follows that $\PP^*$ satisfies all of the axioms for the stratification $X$ and so is quasi-isomorphic to the Deligne sheaf $\mc P^*$. 
\end{proof}

Our second proof is slightly longer but does not require the support and cosupport axioms at all. Rather we  show directly that $\PP^*$ satisfies the axioms Ax1'$(X)$ for any stratification $X$. Consequently, this proof is adaptable to perversities that are not necessarily determined only by codimension; see Remark \ref{R: general}.

\begin{proof}[Proof 2]
We check that $\PP^*$ satisfies Ax1'$(X)$:

\emph{Axiom \ref{I: normal}:} The first two statements are immediate from the definition of $\PP^*$. By construction $\PP^*|_{\mf X-\mf X^{n-2}}\cong \mc E|_{\mf X-\mf X^{n-2}}$ and $\mf X-\mf X^{n-2}\supset X-X^{n-2}$. Finally, $\PP^*$ is $X$-clc  because it is $\mf X$-clc and each stratum of $X$ is contained in a stratum of $\mf X$. 

\emph{Axiom \ref{I: truncate}:} If $x\in S_{n-k}$ then $x$ is contained in a stratum of $\mf X$ of codimension $\ell\leq k$. By construction $H^j(\PP_x)=0$ for $j\geq \bar p(\ell)$. But $\bar p(\ell)\leq \bar p(k)$ as $\bar p$ is a Goresky-MacPherson perversity.

\emph{Axiom \ref{I: attach}':}  
By the growth condition on $\bar p$ and since $\ell\leq k$, $\bar p(k)-\bar p(\ell)\leq k-\ell$, so $n-\ell+\bar p(\ell)+2\geq n-k+\bar p(k)+2=n-\bar q(k)$ and it suffices to demonstrate the vanishing for $j<n-\ell+\bar p(\ell)+2$.

Let $U\cong \R^{n-\ell}\times cL$ be a distinguished neighborhood of $x$ in the $\mf X$ stratification. As all computations are local we may abuse notation, letting  $U= \R^{n-\ell}\times cL$ and letting $\PP^*$ also denote its pullback to this product neighborhood, which remains constructible. Let $\pi_1:\R^{n-\ell}\times cL\to \R^{n-\ell}$ and $\pi_2:\R^{n-\ell}\times cL\to cL$ be the projections, and for some $y\in \R^{n-\ell}$ let $s:cL\into \{y\}\times cL\subset \R^{n-\ell}\times cL$ be the inclusion. By \cite[Proposition 2.7.8]{KS} (letting the $Y_n$ there be close balls in $\R^{n-\ell}$), $\PP^*\cong \pi_2^*R\pi_{2*}\PP^*$. So letting $\mc R_A$ denote the constant $R$ sheaf on the space $A$, we have $$\PP^*\cong\mc R_{\R^{n-\ell}\times cL} \overset{L}{\otimes} \PP^*\cong 
\pi_1^*R_{\R^{n-\ell}}\overset{L}{\otimes} \pi_2^*R\pi_{2*}\PP^*.$$ Now let $x=(y,v)$ with $f_y:\{y\}\into \R^{n-\ell}$ and $f_v:\{v\}\into cL$ the vertex inclusion. By \cite[Remark V.10.20.c]{Bo}, whose hypotheses are satisfied due to the constructibility  \cite[Proposition 4.0.2.2]{Sch03}, $$f^!_{x}\PP^*\cong f^!_y \mc R_{\R^{n-\ell}}\overset{L}{\otimes} f_v^!R\pi_{2*}\PP^*.$$

By \cite[Proposition V.3.7.b]{Bo}, 
$H^i\left(f^!_y \mc R_{\R^{n-\ell}}\right)\cong H^{i-n+\ell}(\mc R_{y})=0$ for $i\leq n-\ell-1$. For $H^i\left(f_v^!R\pi_{2*}\PP^*\right)$, we consider the Cartesian square

\begin{diagram}
\R^{n-\ell}\times cL-\R^{n-\ell}\times \{v\}&\rTo^{\bar \pi_2}&cL-\{v\}\\
\dInto_{\mf i}&&\dInto_{\bar{\mf i}} \\
\R^{n-\ell}\times cL&\rTo^{\pi_2}&cL
\end{diagram}
and the long exact sequence \cite[Section V.1.8]{Bo}
$$\to H^i\left(f_v^!R\pi_{2*}\PP^*\right)\to \H^i\left(cL;R\pi_{2*}\PP^*\right)\xr{\alpha} H^i\left(cL-\{v\}; \bar{\mf i}^*R\pi_{2*}\PP^*\right)\to.$$
We have  $\bar{\mf i}^*R\pi_{2*}\PP^*=R\bar \pi_{2*}\mf i^*\PP^*$; just assume $\PP^*$ injective and look at sections over open sets. So $\alpha$ becomes the restriction $\H^i\left(cL;R\pi_{2*}\PP^*\right)\to H^i\left(cL-\{v\};R\bar \pi_{2*}\mf i^*\PP^*\right)$, which is isomorphic to the attaching map $\H^i(\R^{n-\ell}\times cL;\PP^*)\to  H^i(\R^{n-\ell}\times (cL-\{v\}); \mf i^*\PP^*)$. This is an isomorphism  up through degree $\bar p(\ell)$ by construction. Therefore, $H^i\left(f_v^!R\pi_{2*}\PP^*\right)=0$ for $i\leq\bar p(\ell)$. It is also $0$ in degree $\bar p(\ell)+1$ as $\H^{\bar p(\ell)+1}\left(cL;R\pi_{2*}\PP^*\right)\cong \H^{\bar p(\ell)+1}\left(\R^{n-\ell}\times cL;\PP^*\right)\cong H^{\bar p(\ell)+1}(\PP^*_x)=0$; the last equality by construction and the middle isomorphism by the constructibility and \cite[Proposition 4.0.2.2]{Sch03}.
So by the K\"unneth Theorem, 
$H^j\left(f_x^!\PP^*\right)=0$ for $j<n-\ell+\bar p(\ell)+2$ as desired.
\end{proof}

\begin{remark}\label{R: twisted}
Both of our proofs can be generalized to non-constant coefficient systems by using instead maximally coarse filtrations that depend on the coefficients. See \cite[Section 3]{HS91} and \cite{GBF43} for details. 
\end{remark}

\begin{remark}\label{R: general}
It is also possible to characterize axiomatically Deligne sheaves with arbitrary perversities $\bar p:\{\text{singular strata}\}\to \Z$ that in particular do not necessarily depend only on codimension \cite{GBF23}.  Our second proof can be applied more generally to show that if $\mc X$ is a CS set coarsening $X$ with respective perversities $\bar p_{\mc X}$ and $\bar p_X$, then the respective Deligne sheaves $\mc P^*_{X,\bar p_X}$ and $\mc P^*_{\mc X,\bar p_{\mc X}}$
are quasi-isomorphic if whenever $S$ is a stratum of $X$ contained in the stratum $\mc S$ of $\mc X$ then we have  $\bar p_{\mc X}(\mc S)\leq  \bar p_{X}(S)\leq\bar p_{\mc X}(\mc S)+\codim(S)-\codim(\mc S)$ if $\mc S$ is singular and $0\leq \bar p_X(S)\leq \codim(S)-2$ if $\mc S$ is regular. 
This generalizes recent results of Chataur, Saralegi-Aranguren, and Tanre  \cite{CST-inv}. Again, see \cite{GBF43} for details.
\end{remark}

\bibliographystyle{amsplain}
\providecommand{\bysame}{\leavevmode\hbox to3em{\hrulefill}\thinspace}
\providecommand{\MR}{\relax\ifhmode\unskip\space\fi MR }
\providecommand{\MRhref}[2]{%
  \href{http://www.ams.org/mathscinet-getitem?mr=#1}{#2}
}
\providecommand{\href}[2]{#2}

\end{document}